\documentclass[preprint,11pt,authoryear]{elsarticle}
\usepackage{amssymb}
 \usepackage{amsmath}
 \usepackage[latin1]{inputenc}
 \usepackage[T1]{fontenc}
\usepackage{times}
 \usepackage{dsfont}
 \usepackage{natbib}
 \usepackage{mathrsfs}
 \usepackage{hyperref}
 \usepackage[english]{babel}
 \input{bb.sty}
\date{}

 \usepackage{graphicx}

\newtheorem{theorem}{Theorem}[section]
\newtheorem{lemma}[theorem]{Lemma}

\newtheorem{e-proposition}[theorem]{Proposition}

\newtheorem{definition}[theorem]{Definition\rm}

\newtheorem{theoreme}{Th\'eor\`eme}[section]

\newtheorem{corollaire}[theoreme]{Corollaire}
\newtheorem{e-definition}[theoreme]{D\'efinition\rm}

\renewcommand{\theequation}{\arabic{equation}}
\def \ds{\displaystyle}
\begin{document}

\begin{frontmatter}

%

\title{Nonparametric regression estimation for quasi-associated
Hilbertian processes}

 \author[]{Lahcen DOUGE}
 \address{LaMSAFA, FSTG, University Cadi Ayyad, B.P. 549 Marrakech Morocco}
 
\begin{abstract}
 We establish  the asymptotic normality of the kernel type estimator for the regression function constructed from
 quasi-associated data when the explanatory variable takes
 its values in a separable Hilbert space.
\end{abstract}


\end{frontmatter}


\section{Introduction} \noindent The study of statistical models adapted to infinite dimensional data has been the subject of several works in
the recent statistical literature (see \cite{Bosq},  \cite{Ramsay},  \cite{Ferraty}.\\
In this paper we investigate nonparametric estimation of the regression function  when the explanatory variable is
functional and taking values in a separable Hilbert space  and the response is scalar. We establish the asymptotic normality of the Nadaraya-Watson type estimator for the regression functional for
quasi-associated processes. The asymptotic properties of this
estimator have been studied by \cite{Ferraty} and
\cite{Masry} in the case of strongly mixing processes.\\
The concept of quasi-association was introduced for real-valued random
fields by \cite{Bulinski} and provides a unified approach to studying both families of positively or negatively associated and gaussian random variables. This notion is a special case of weak dependence introduced by \cite{Doukhan}  for real-valued random variables.\\
Before recalling the definition of quasi-association for real random vectors, denote by
 \begin{eqnarray*}
 	\mathrm{Lip}(h) =\sup_{x\neq y}\frac{|h(x)-h(y)|}{\|x-y\|_{\mathcal{F}}}
 \end{eqnarray*}
 the Lipschitz modulus of a function $h:\mathcal{F} \rightarrow \mathds{R}$, where $\mathcal{F}$ is some normed space with norm $\|.\|_{\mathcal{F}}$.
\begin{definition}\label{4def1}
A sequence of random vectors $(X_i)_{i\in \mathds{N}}$ with values
in $\mathds{R}^{d}$, $d\geq 1$, is said to be quasi-associated if for any finite and 
disjoint subsets $I, J\subset \mathds{N}$ and all bounded Lipschitz
functions $ f:\mathds{R}^{|I|d}\rightarrow
\mathds{R},~g:\mathds{R}^{|J|d}\rightarrow \mathds{R}$, one has
\begin{eqnarray*}
\Big|\mathrm{Cov}\big(f(X_i,i\in I),g(X_j, j\in
J)\big)\Big|\leq \mathrm{Lip}(f)\mathrm{Lip}(g)\sum_{i\in
I}\sum_{j\in
J}\sum_{k,l=1}^{d}\big|\mathrm{Cov}(X_{i}^{k},X_{j}^{l})\big|,
\end{eqnarray*}
where $|I|$ denotes the cardinality of $I$, $X_{i}^{k}$ denotes the kth component
of $X_i$.
\end{definition}
Now we introduce a definition of quasi-association for random variables with values in a separable Hilbert space similar to the notion of weakly association for Hilbertian processes in \cite{Burton}.
\begin{definition}
A sequence $(Z_i)_{i\in \mathds{N}}$ of random variables taking
values in a separable Hilbert space $\big(\mathcal{E},<.,.>\big)$ is called
quasi-associated, with respect to an orthonormal basis
$\{e_k,\,k\geq 1\}$ in $\mathcal{E}$, if for any $d\geq 1$, the
$d$-dimensional sequence $\Big\{\big(<Z_i,e_1>,\ldots,<Z_i,e_d>\big),\,i\in
\mathds{N}\Big\}$ is quasi-associated.
\end{definition}
Some examples of quasi-associated Hilbertian processes are given in \cite{Douge}.\\ 
Let $\mathcal{H}$ be a separable real Hilbert space with
the norm $\|\cdot\|$ generated by an inner product
$<\cdot,\cdot>$. Let $(X,Y),\,(X_1,Y_1),\,(X_2,Y_2),\ldots $ be a sequence of stationary quasi-associated and identically distributed random variables in the separable Hilbert space $\mathcal{E}:=\mathcal{H}\times \mathds{R}$ and let $x$ be a fixed element of $\mathcal{H}$. the regression operator of $\varphi(Y)$ on $X$ is defined by
\begin{eqnarray*}
r(x):=\mathrm{E}\big(\varphi(Y)|X=x\big),
\end{eqnarray*}
where $\varphi$ is a real-valued Borel function defined on
$\mathds{R}$. The estimate we consider here is of kernel type defined by
\begin{eqnarray*}
r_n(x)=\frac{\sum_{i=1}^{n}\varphi(Y_i)
K\left(\frac{\|x-X_i\|}{h_n}\right)}{\sum_{i=1}^{n}
K\left(\frac{\|x-X_i\|}{h_n}\right)},
\end{eqnarray*}
where $K$ is a kernel function and $\{h_n\}_{n\geq 0}$ is a sequence of positive constants such that, as $n\rightarrow \infty$, $h_n \rightarrow 0$ and $nh_n \rightarrow \infty$. \\
For $i=1,\ldots,n$. Set
\begin{eqnarray*}
\Delta_i(x):=K\left(\frac{\|x-X_i\|}{h_n}\right) \quad\text{and}\quad
\Gamma_i(x):= \varphi(Y_i) K\left(\frac{\|x-X_i\|}{h_n}\right).
\end{eqnarray*}
Now define
\begin{eqnarray*}
f_n(x):=\frac{1}{n \mathrm{E}
\Delta_1(x)}\sum_{i=1}^{n}\Delta_i(x)\quad\text{and}\quad g_n(x):=\frac{1}{n
\mathrm{E}\Delta_1(x)}\sum_{i=1}^{n}\Gamma_i(x).
\end{eqnarray*}
So that $r_n(x)=g_n(x)/f_n(x)$.\\
Define also the truncated kernel estimator of $r(x)$ by
$\hat{r}_n(x)=\hat{g}_n(x)/f_n(x)$, where
\begin{eqnarray*}
\hat{g}_n(x):=\frac{1}{n\mathrm{E}\Delta_1(x)}\sum_{i=1}^n
\hat{\Gamma}_i(x)\quad\text{and}\quad
\hat{\Gamma}_i(x):=\varphi(Y_i)\mathds{1}_{\{|\varphi(Y_i)|\leq
b_n\}} \Delta_i(x),
\end{eqnarray*}
where $\{b_n\}_{n\geq 0}$ denoting a positive sequence such that
$b_n\rightarrow \infty$. \\
 Denote
\begin{eqnarray*}
\lambda_{ij}:=\sum_{k=1}^{\infty}\sum_{l=1}^{\infty}\big|\mathrm{Cov}(X_{i}^{k},X_{j}^{l})\big|+
\sum_{k=1}^{\infty}\big|\mathrm{Cov}(X_{i}^{k},Y_j)\big|+\sum_{l=1}^{\infty}\big|\mathrm{Cov}(Y_i,X_{j}^{l})\big|+\big|\mathrm{Cov}(Y_i,Y_j)\big|,
\end{eqnarray*}
where $X_{i}^{k}:=<X_i,e_k>$ and $\ds \lambda_k=\sup_{s\geq
k}\sum_{|i-j|\geq s}\lambda_{i,j}$.\\
Let $D(x):=\|x-X_1\|$ a
real-valued nonnegative random variable. Denote its
distribution by $F(u,x):=P\big(D(x)\leq u\big)$, $u\in\mathds{R}$. \cite{Gasser} assume that  if there exist a function $\phi$ such that $\lim_{u\rightarrow 0} F(u,x)/\phi(u) =:f(x)$, then they refer to $f$ as a probability density
of $X_1$.
\section{Assumptions and main results}
\noindent As usual in nonparametric functional estimation
problems, we introduce the following assumptions which we need to establish the main result.\\\\
\noindent \textbf{Assumptions}
\begin{itemize}
\item[A1 (i)]There exist some constants $c_1,\, c_2,\, c_3,\, c_4 >0$ such that for $0<u<c_1$
\begin{eqnarray*}
0<c_2 \phi(u) f_1(x)\leq F(u,x)\leq c_3 \phi(u)
f_1(x),
\end{eqnarray*}
where $\phi(u)\rightarrow 0$ as $u\rightarrow 0$ and $f_1$ is a
function from $\mathcal{H}$ to $\mathds{R}^{+}$. 
 \item[(ii)]
\begin{displaymath}
\sup_{i\neq j}P\big(D_i(x)\leq u,\, D_j(x)\leq u\big)\leq
c_4\psi (u)f_2(x),
\end{displaymath}
 where $\psi(u)\rightarrow 0$ as $u\rightarrow
0$ and $f_2$ is a function from $\mathcal{H}$ to $\mathds{R}^{+}$. We
assume that the ratio $\psi(u)/\phi(u)^2$ is bounded.
\item[A2 (i)] $0<c_5\mathds{1}_{[0,1]}\leq K\leq
c_6\mathds{1}_{[0,1]}$ for some constants $c_5$ and $c_6$. \item[(ii)] $K$ is a Lipschitz function. 
\item[(iii)] $\phi$ is derivable and $\lim_{u\rightarrow 0} \ds
    \frac{u}{\phi(u)}\int_{0}^{1}K^{j}(y)\phi^{'}(uy)\mathrm{d}y:=
    C_j$, $j=1,2$.
\item[A3 (i)]$\varphi$ is a Lipschitz function.
\item[(ii)]  $\mathrm{E}\exp\big(|\varphi(Y_1)|\big)<\infty$.
\item[(iii)]$\big|r(u)-r(v)\big|\leq c_7 \|u-v\|^{\beta}$, $u,v
    \in \mathcal{H}$, for some $\beta>0$ and $c_7>0$.
\item[(iv)] The function $g_2(u):=
    \mathrm{E}\big(\varphi^2(Y_1)|X_1=u\big),~u\in \mathcal{H}$, exists
    and is uniformly continuous in some
    neighborhood of $x$.
\end{itemize}
\noindent Assumption A1(i) is inspired from the work of \cite{Gasser} and assumption A1(ii) gives the behavior of the joint distribution $\big(D_i(x), D_j(x)\big)$.
Assumptions A2(i)(ii) are standard for $K$ and assumption A2(iii) is necessary to obtain an expression of the asymptotic variance. Assumptions A3(i)(iii)(iv) are a mild smoothness assumptions on the regression functional $r$ and the function $\varphi$.\\

\begin{theorem}\label{prop4.3}
Suppose that conditions  (A1)-(A3)  hold. Suppose in addition that $\phi(h_n)=\mathcal{O}(h_n^b)$,  $(\log n)^2\phi(h_n)^{a\delta-(1+\frac{2}{b})}\rightarrow 0$ and $n\phi(h_n)^{1+2\delta}\rightarrow \infty$ for some $0<\delta<1$ and $b>0$. If $\lambda_k=\mathcal{O}(k^{-a})$ for some $a>\ds\frac{2+b}{\delta b}$, then
\begin{eqnarray}\label{eqq4.1}
\sqrt{n\phi(h_n)}\big(g_n(x)-\mathrm{E}g_n(x)\big)\xrightarrow{d}
\mathcal{N}\left(0,\sigma_1^2(x)\right)
\end{eqnarray}
 and
\begin{eqnarray*}
\sqrt{n\phi(h_n)}\Big(\big[g_n(x)-r(x)f_n(x)\big]-\mathrm{E}\big[g_n(x)-r(x)f_n(x)\big]\Big)\xrightarrow{d}
\mathcal{N}\left(0,\sigma_2^2(x)\right)\quad \mbox{as}\,\, n\rightarrow \infty,
\end{eqnarray*}
where $\displaystyle
\sigma_1^2(x):=\frac{C_2}{C_1^2}\frac{g_2(x)}{f_1(x)}$\quad and
\quad $\displaystyle
\sigma_2^2(x):=\frac{C_2}{C_1^2}\frac{g_2(x)-r^2(x)}{f_1(x)}\cdot$
\end{theorem}
\vspace{0.3cm}

\begin{corollaire}\label{thm4.2}
Suppose that conditions  (A1)-(A3) are satisfied. Suppose in addition that $\phi(h_n)=\mathcal{O}(h_n^b)$ and $nh_n^{2\beta}\phi(h_n)\rightarrow 0$,  $n\phi(h_n)^{1+2\delta}\rightarrow \infty$ for some $0<\delta<1$ and $b>0$. If $\lambda_k=\mathcal{O}(k^{-a})$ for some
$a>\ds\frac{2+b}{\delta b}$, then
\begin{eqnarray*}
\sqrt{n\phi(h_n)}\big(r_n(x)-r(x)\big)\xrightarrow{d}
\mathcal{N}\big(0,\sigma_2^2(x)\big)\quad \mbox{as}\,\, n\rightarrow \infty.
\end{eqnarray*}
\end{corollaire}
\section{Proofs}
\noindent In order to prove Theorem \ref{prop4.3}, we first introduce the following two lemmas. Denote by $BL(\mathcal{E}^m)$, with $m$ a strictly
positive integer, the set of Lipschitz and bounded functions
$f:\mathcal{E}^m \rightarrow \mathds{R}$. We equip $\mathcal{E}^m$ with the norm $\|x\|_{\mathcal{E}^{m}}=\sum_{s=1}^{m}\|x_s\|_{\mathcal{E}}$, where $x=(x_1,\ldots,x_m)\in \mathcal{E}^{m}$ and $\|.\|_{\mathcal{E}}$ is the norm induced by the inner product on $\mathcal{E}$. Throughout the demonstrations, denote by $C$ different constants whose values are allowed to change. 
\begin{lemma}\label{lem4.1}
\noindent Let $(Z_n)_{ n\in \mathds{N}}$ be a quasi-associated
sequence of random variables with values in $\mathcal{E}$. Let
$f\in BL\big(\mathcal{E}^{|I|}\big)$, $g\in
BL\big(\mathcal{E}^{|J|}\big)$, for some finite
disjoint subsets $I,~ J \subset\mathds{N}$. Then
\begin{eqnarray*}
\Big|\mathrm{Cov}\big(f(Z_i,i\in I),\,g(Z_j,j\in
J)\big)\Big|\leq \mathrm{Lip}(f) \mathrm{Lip}(g)
\sum_{i\in I} \sum_{j\in J} \sum_{k=1}^{\infty}\sum_{l=1}^{\infty}
\big|\mathrm{Cov}(Z_{i}^{k},Z_{j}^{l})\big|.
\end{eqnarray*}
\end{lemma}
\noindent \textbf{Proof.} Let$\{e_k,~k\geq 1\}$
is an orthonormal basis in $\mathcal{E}$ and let $F:\mathds{R}^{n|I|}\rightarrow \mathds{R}$ and
$G:\mathds{R}^{n|J|}\rightarrow \mathds{R}$ be two functions such
that
\begin{eqnarray*}
 F\big(<x_i,e_k>,~
 1\leq k\leq n,~ i\in I\big)=f\left(\sum_{k=1}^n <x_i,e_k> e_k,~i\in I\right),\quad x_i\in \mathcal{E},\, i\in I 
\end{eqnarray*}
 and
\begin{eqnarray*}
 G\big(<x_i,e_k>,~ 1\leq k\leq n,~ j\in J\big)=g\left(\sum_{k=1}^n <x_i,e_k> e_k,~j\in
J\right),\quad x_j\in \mathcal{E},\, j\in J. 
\end{eqnarray*}
$f$ and $g$ are continuous and bounded, then by the
Lebesgue's dominated convergence theorem
\begin{eqnarray}\label{eq4.2}
\nonumber\lefteqn{\left|\mathrm{Cov}\left(f\left(X_i,i\in
I\right),g\left(X_j,j\in
J\right)\right)\right|}\\&=&\nonumber\lim_{n\rightarrow \infty}
\left|\mathrm{Cov}\left(f\left(\sum_{k=1}^n X_{i}^{k} e_k,~i\in
I\right),g\left(\sum_{k=1}^n X_{j}^{k} e_k,~j\in
J\right)\right)\right|\\
&=&\lim_{n\rightarrow \infty} \Big|\mathrm{Cov}\big(F(X_{i}^{k},~
1\leq k\leq n,~ i\in I),G(X_{j}^{k},~ 1\leq k\leq n,~ j\in
J)\big)\Big|.
\end{eqnarray}
Next, by quasi-association of the sequence $(X_n)$, we have
\begin{eqnarray}\label{eq4.3}
\nonumber\lefteqn{\Big|\mathrm{Cov}\big(F(X_{i}^{k},~ 1\leq k\leq
n,~ i\in I),G(X_{j}^{k},~ 1\leq k\leq n,~ j\in
J)\big)\Big|}\\&\leq& \nonumber\mathrm{Lip}(F)\mathrm{Lip}(G)
\sum_{i\in I}\sum_{j\in J}\sum_{k=1}^n\sum_{l=1}^n
|\mathrm{Cov}(X_{i}^{k},X_{j}^{l})|\\&\leq&
\mathrm{Lip}(f)\mathrm{Lip}(g) \sum_{i\in I}\sum_{j\in
J}\sum_{k=1}^n\sum_{l=1}^n |\mathrm{Cov}(X_{i}^{k},X_{j}^{l})|.
\end{eqnarray}
The proof is completed by (\ref{eq4.2}) and
(\ref{eq4.3}).\cqfd

\begin{lemma}\label{prop4.1}
Suppose that conditions  (A1)-(A3)  hold. Suppose  that $\lambda_k=\mathcal{O}(k^{-a})$ for some $a>1+\dfrac{2}{b}$. If $\phi(h_n)=\mathcal{O}(h_n^b)$ and $ (\log n)^{4/\gamma} \phi(h_n) \rightarrow 0$,  for $\gamma =1-\dfrac{2+b}{ab}$,  as $n\rightarrow \infty$
. Then
\begin{eqnarray*}
n\phi(h_n)\mathrm{Var}(g_n(x))=\sigma_1^2(x)+o(1)
\end{eqnarray*}
and
\begin{eqnarray*}
n\phi(h_n)\mathrm{Var}(g_n(x)-r(x)f_n(x))=\sigma_2^2(x)+
o(1).
\end{eqnarray*}
\end{lemma}

\noindent\textbf{Proof.}
\begin{align*}
\mathrm{Var}g_n(x)=&~\mathrm{E}\left[g_n(x)-\hat{g}_n(x)-\mathrm{E}\left(g_n(x)-\hat{g}_n(x)\right)\right]^2 + \mathrm{Var}\hat{g}_n(x)\\
&+
2\mathrm{E}\left[\left(g_n(x)-\hat{g}_n(x)-\mathrm{E}\left(g_n(x)-\hat{g}_n(x)\right)\right)\left(\hat{g}_n(x)-\mathrm{E}\hat{g}_n(x)\right)
\right].
\end{align*}
By A3 (ii) and by using H\"{o}lder inequality, it follows
that
\begin{eqnarray}\label{eq4.51}
\lefteqn{\left(n \phi(h_n)\right)\mathrm{E}\left|g_n(x)-\hat{g}_n(x)-\mathrm{E}\left(g_n(x)-\hat{g}_n(x)\right)\right|^2}\nonumber\\
 &\leq&
C^2\frac{n^2}{n\phi(h_n)}\mathrm{E}\left[|\varphi(Y_1)|^2\mathds{1}_{\{|\varphi(Y_1)|>b_n\}}
\Delta_1^2(x)\right]
\nonumber\\
&\leq& C^q\frac{n^2}{n
\phi(h_n)}\left(\mathrm{E}|\varphi(Y_1)|^{4}\right)^{1/2}
n^{-b_0/2}
\left(\mathrm{E}\left(\exp\left|\varphi(Y_1)\right|\right)\right)^{1/2}\rightarrow
0,
\end{eqnarray}
for $b_n=b_0 \log(n)$, $b_0>0$, large enough. 

\begin{eqnarray}\label{eq4.9}
\mathrm{Var}\hat{g}_n(x)=\frac{1}{n\left[\mathrm{E}\Delta_1(x)\right]^2}\mathrm{Var}\hat{\Gamma}_1(x)+
\frac{1}{n^2\left[\mathrm{E}\Delta_1(x)\right]^2}\sum_{i=1}^n\sum_{\substack{j=1\\j\neq
i}}^{n}\mathrm{Cov}\left(\hat{\Gamma}_i(x),\hat{\Gamma}_j(x)\right).
\end{eqnarray}

\begin{align*}
\mathrm{Var}\hat{\Gamma}_1(x)=&~\mathrm{E}\left(\varphi^2(Y_1)\mathds{1}_{\{|\varphi(Y_1)|\leq
b_n\}} \Delta_1^2(x)\right)-
\left(\mathrm{E}\left(\varphi(Y_1)\mathds{1}_{\{|\varphi(Y_1)|\leq b_n\}} \Delta_1(x)\right)\right)^2\\
=&~\mathrm{E}\left(\varphi^2(Y_1) \Delta_1^2(x)\right)- \mathrm{E}\left(\varphi^2(Y_1)\mathds{1}_{\{|\varphi(Y_1)|>b_n\}}\Delta_1^2(x)\right)\\
&- \left(\mathrm{E}\left(\varphi(Y_1)
\Delta_1(x)\right)-\mathrm{E}\left(\varphi(Y_1)\mathds{1}_{\{|\varphi(Y_1)|>b_n\}}
\Delta_1(x)\right)\right)^2.
\end{align*}
Now,
\begin{eqnarray*}
\mathrm{E}\left(\varphi^2(Y_1)
\Delta_1^2(x)\right)&=&\mathrm{E}\left(g_2(X_1)\Delta_1^2(x)\right)\\
&=&g_2(x)\mathrm{E}\left(\Delta_1^2(x)\right)+\mathrm{E}\left(\left(g_2(X_1)-g_2(x)\right)\Delta_1^2(x)\right).
\end{eqnarray*}
 By condition A3(iv),
\begin{eqnarray*}
\mathrm{E}\left(\left(g_2(X_1)-g_2(x)\right)\Delta_1^2(x)\right) &\leq& \sup_{\{u:\|x-u\|\leq h\}}|g_2(u)-g_2(x)|\mathrm{E}( \Delta_1^2(x))\\
&=& o(1) \mathrm{E}\Delta_1^2(x).
\end{eqnarray*}
Thus,
\begin{eqnarray*}
\mathrm{E}\left(\varphi^2(Y_1)\Delta_1^2(x)\right)=(g_2(x)+o(1))\mathrm{E}\Delta_1^2(x).
\end{eqnarray*}
By condition A3(iii),
\begin{eqnarray*}
\mathrm{E}\left(\varphi(Y_1)
\Delta_1(x)\right)&=&r(x)\mathrm{E}\Delta_1(x)+\mathrm{E}\left(\left(r(X_1)-r(x)\right)\Delta_1(x)\right)\\
&=& (r(x)+\mathcal{O}(h_n^{\beta}))\mathrm{E}\Delta_1(x).
\end{eqnarray*}
For $j=1,2$,
\begin{eqnarray*}
\frac{1}{\phi(h_n)}\mathrm{E}[\Delta_1^j(x)]&=&\frac{1}{\phi(h_n)}\int_{0}^{h_n}K^{j}(u/h_n)\mathrm{d}F(u;x)
\rightarrow C_jf_1(x).
\end{eqnarray*}
On the other hand, using H\"{o}lder inequality, we have
\begin{eqnarray*}
\mathrm{E}\left(\varphi^2(Y_1)\mathds{1}_{\{|\varphi(Y_1)|>b_n\}}\Delta_1^2(x)\right)\leq
C b_n^{2-s}
\end{eqnarray*}
and
\begin{eqnarray*}
\left|\mathrm{E}\left(\varphi(Y_1)\mathds{1}_{\{|\varphi(Y_1)|>b_n\}}\Delta_1(x)\right)\right|\leq
C b_n^{1-s}.
\end{eqnarray*}
Then,
\begin{eqnarray}\label{eq4.10}
\frac{\phi(h_n)}{\left[\mathrm{E}\Delta_1(x)\right]^2}\mathrm{Var}\hat{\Gamma}_1(x)\rightarrow
\frac{C_2}{C_1^2}\frac{g_2(x)}{f_1(x)}.
\end{eqnarray}
\vspace{0.3cm}

\noindent Next, we start by studying the sum in the right hand
side of (\ref{eq4.9}). We will use the following natural
decomposition, in which $(v_n)$ is some sequence of positive
integers
\begin{eqnarray}\label{equa4.3}
\sum_{i=1}^n\sum_{\substack{j=1\\j\neq
i}}^{n}\mathrm{Cov}\left(\hat{\Gamma}_i(x),\hat{\Gamma}_j(x)\right)=\sum_{i=1}^n\sum_{\substack{j=1\\0<|i-j|\leq
v_n}}^{n}\mathrm{Cov}\left(\hat{\Gamma}_i(x),\hat{\Gamma}_j(x)\right)+\sum_{i=1}^n\sum_{\substack{j=1\\|i-j|>v_n}}^{n}
\mathrm{Cov}\left(\hat{\Gamma}_i(x),\hat{\Gamma}_j(x)\right).
\end{eqnarray}
 It follows from the assumptions A1, A2(i), for $i\neq j$, that
\begin{eqnarray*}\label{eq4.5}
\left|\mathrm{E}\left(\hat{\Gamma}_i(x) \hat{\Gamma}_j(x)\right)\right|&\leq& b_n^2\mathrm{E}\left| \Delta_i(x)\Delta_j(x)\right|\nonumber\\
&\leq& C b_n^2\sup_{i\neq j}P[(X_i,X_j)\in \mathcal{B}(x,h_n)\times\mathcal{B}(x,h_n)]\nonumber\\
&\leq& C b_n^2\phi(h_n)^2
\end{eqnarray*}
and
\begin{eqnarray*}\label{eq4.6}
\left|\mathrm{E}\left(\hat{\Gamma}_1(x)\right)\right|&\leq& b_n\mathrm{E}\left|\Delta_1(x)\right|\nonumber\\
&\leq& b_n P[X_1\in \mathcal{B}(x,h_n)]\nonumber\\
&\leq& C b_n\phi(h_n).
\end{eqnarray*}
Then
\begin{eqnarray*}
\nonumber\sum_{i=1}^n\sum_{\substack{j=1\\0<|i-j|\leq
v_n}}^{n}\left|\mathrm{Cov}\left(\hat{\Gamma}_i(x),\hat{\Gamma}_j(x)\right)\right|
&\leq&nv_n \left[ \max_{i\neq j}\left|\mathrm{E}\left(\hat{\Gamma}_i(x) \hat{\Gamma}_j(x)\right)\right|+\left(\mathrm{E}\hat{\Gamma}_1(x)\right)^2\right]\\
&\leq& C n v_n b_n^2\phi(h_n)^2.
\end{eqnarray*}
By applying Lemma \ref{lem4.1} to the sequence $(Z_i)$, we get
\begin{eqnarray*}
\nonumber\sum_{i=1}^n\sum_{\substack{j=1\\|i-j|>v_n}}^{n}
\left|\mathrm{Cov}\left(\hat{\Gamma}_i(x),\hat{\Gamma}_j(x)\right)\right|&\leq&
C b_n^2 h_n^{-2}
\sum_{i=1}^n\sum_{\substack{j=1\\|i-j|>v_n}}^{n}\lambda_{ij}\\\nonumber
&\leq& C n b_n^2 h_n^{-2}\lambda_{v_n}.
\end{eqnarray*}
By this choice of $(b_n)$, we obtain
\begin{eqnarray*}
\frac{1}{n\phi(h_n)}\sum_{i=1}^n\sum_{\substack{j=1\\i\neq j}}^{n}
\mathrm{Cov}\left(\hat{\Gamma}_i(x),\hat{\Gamma}_j(x)\right)&\leq& C\left[b_n^2\phi(h_n)v_n+b_n^2\phi(h_n)^{-1}h_n^{-2}v_n^{-a}\right]\\
&\leq&
C\left[\log^2(n)\phi(h_n)v_n+\log^2(n)\phi(h_n)^{-(1+\frac{2}{b})}v_n^{-a}\right],
\end{eqnarray*}
Let $\kappa=1-\dfrac{\gamma}{2}$, then by choosing $v_n$ such that $v_n=[\phi(h_n)^{-\kappa}]$, we deduce
that
\begin{eqnarray}\label{eq4.12}
\frac{1}{n\phi(h_n)}\sum_{i=1}^n\sum_{\substack{j=1\\i\neq j}}^{n}
\mathrm{Cov}\left(\hat{\Gamma}_i(x),\hat{\Gamma}_j(x)\right)\rightarrow
0.
\end{eqnarray}


Then,
\begin{eqnarray}\label{equat4.1}
n\phi(h_n)\mathrm{Var}\hat{g}_n(x)\rightarrow
\frac{C_2}{C_1^2}\frac{g_2(x)}{f_1(x)}\cdot
\end{eqnarray}
The proof is completed by showing that
\begin{eqnarray*}
\lefteqn{(n\phi(h_n))\mathrm{E}\left[\left(g_n(x)-\hat{g}_n(x)-\mathrm{E}\left(g_n(x)-\hat{g}_n(x)\right)\right)\left(\hat{g}_n(x)-\mathrm{E}\hat{g}_n(x)\right) \right]}\\
&\leq&\left[(n\phi(h_n))\mathrm{E}\left(g_n(x)-\hat{g}_n(x)-\mathrm{E}\left(g_n(x)-\hat{g}_n(x)\right)\right)^2\right]^{1/2}
\left[\left(n\phi(h_n)\right)\mathrm{Var}\hat{g}_n(x)\right]^{1/2}\rightarrow
0.
\end{eqnarray*}
The second assertion follows by replacing $\varphi(Y_i)$ with
$\varphi(Y_i)-r(x)$.\cqfd \vspace{0.5cm}


\noindent\textbf{Proof of Theorem \ref{prop4.3}.} The basic technique in establishing (\ref{eqq4.1})
consists in spliting the set $\{1,\ldots, n\}$ into $k$ large
$p$-blocks and small $q$-blocks, to be denoted by $I_j$ and $J_j$,
$j=1,\ldots,k$, respectively as follows:
\begin{eqnarray*}
I_j&=&\{(j-1)(p+q)+1,\ldots,(j-1)(p+q)+p\},\\
J_j&=&\{(j-1)(p+q)+p+1,\ldots,j(p+q)\},
\end{eqnarray*}
where $p=p_n,~q=q_n$ are positive integers tending to $\infty$, as
$n\rightarrow \infty$, and $k=k_n$ is defined by $k=[n/(p+q)]$ ($[x]$ stands for the integral part of $x$).\\
We suppose that
\begin{eqnarray*}
\frac{qk}{n}\rightarrow 0 \quad \mbox{and} \quad
\frac{pk}{n}\rightarrow 1.
\end{eqnarray*}
Define
\begin{eqnarray*}
Z_{ni}(x)=\frac{\phi(h_n)^{1/2}}{\sqrt{n}~\mathrm{E}\left(\Delta_1(x)\right)}\left(\hat{\Gamma}_i(x)-\mathrm{E}\hat{\Gamma}_i
(x)\right).
\end{eqnarray*}
and
\begin{eqnarray*}
S_n:=\sum_{i=1}^n Z_{ni}(x).
\end{eqnarray*}
Now, we want to show that
\begin{eqnarray}\label{eq4.13}
\sqrt{n\phi(h_n)}\left(g_n(x)-\mathrm{E}g_n(x)\right)-S_n\xrightarrow{\mathbb{L}^2}
0
\end{eqnarray}
and
\begin{eqnarray}\label{eq4.14}
S_n\xrightarrow{\mathcal{D}} \mathcal{N}(0,\sigma_1^2(x)).
\end{eqnarray}
(\ref{eq4.13}) is proved as in (\ref{eq4.51}) by using H\"{o}lder inequality and by choosing $b_n=b_0\log(n)$.\\
 Let us prove (\ref{eq4.14}).\\
For $j=1,\ldots,k$, let $\eta_j$, $\xi_j$, $\zeta_k$ be defined as
follows
\begin{eqnarray*}
\eta_j:=\sum_{i\in I_j}Z_{ni}(x),\quad \xi_j:=\sum_{i\in
J_j}Z_{ni}(x),\quad \zeta_k:=\sum_{i=k(p+q)+1}^n Z_{ni}(x),
\end{eqnarray*}
so that
\begin{eqnarray*}
S_n=\sum_{j=1}^k \eta_j+\sum_{j=1}^k \xi_j+\zeta_k=:
T_n+T_n^{'}+T_n^{''}.
\end{eqnarray*}
Convergence (\ref{eq4.14}) will be established by showing that
\begin{eqnarray}\label{eq4.15}
T_n \xrightarrow{\mathcal{D}}\mathcal{N}(0,\sigma_1^2(x))
\end{eqnarray}
and
\begin{eqnarray}\label{eq4.16}
\mathrm{E}(T_n^{'})^2 + \mathrm{E}(T_n^{''})^2\rightarrow 0.
\end{eqnarray}
The proof of convergence in (\ref{eq4.15}) consists in using
characteristic functions and showing the following two results
\begin{eqnarray}\label{eq4.17}
\Big|\mathrm{E}\left(e^{it\sum_{j=1}^k
\eta_j}\right)-\prod_{j=1}^{k}\mathrm{E}\left(e^{it
\eta_j}\right)\Big|\rightarrow 0,
\end{eqnarray}
and
\begin{eqnarray}\label{eq4.18}
k\mathrm{Var}(\eta_1)\rightarrow \sigma_1^2(x),\quad
k\mathrm{E}\left(\eta_1^{2}\mathds{1}_{\{\eta_1^{}>\varepsilon
\sigma_1(x) \}}\right)\rightarrow 0.
\end{eqnarray}
(\ref{eq4.18}) is the standard Lindeberg-Feller condition for
asymptotic normality of $T_n$ under independence. \vspace{0.3cm}

\noindent \emph{Proof of convergence in (\ref{eq4.16})}. By using stationarity,
\begin{eqnarray}\label{eq4.19}
\mathrm{E}[T_n^{'}]^2\leq k\mathrm{Var}(\xi_1)+2\sum_{1\leq
i<j\leq k} |\mathrm{Cov}(\xi_i,\xi_j)|
\end{eqnarray}
and
\begin{eqnarray}\label{eq4.20}
k\mathrm{Var}(\xi_1)\leq qk\mathrm{Var}(Z_{n1}(x))+2k\sum_{1\leq
i<j\leq q} |\mathrm{Cov}(Z_{ni}(x),Z_{nj}(x))|.
\end{eqnarray}
The first term in the right-hand side of (\ref{eq4.20}) can be
treated by means of (\ref{eq4.10})
\begin{eqnarray}\label{eq4.21}
qk\mathrm{Var}Z_{n1}(x)=\frac{qk\phi(h_n)}{n\left[\mathrm{E}\Delta_1(x)\right]^2}\mathrm{Var}\hat{\Gamma}_1(x)\rightarrow
0,\quad n\rightarrow \infty.
\end{eqnarray}
The second term can be treated by using the same decomposition
given in (\ref{equa4.3})
\begin{eqnarray*}
k\sum_{1\leq i<j\leq q}
|\mathrm{Cov}(Z_{ni},Z_{nj})|&=&\frac{k\phi(h_n)}{n\left[\mathrm{E}\Delta_1(x)\right]^2}
\sum_{1\leq i<j\leq q}\left|\mathrm{Cov}\left(\hat{\Gamma}_i(x),\hat{\Gamma}_j(x)\right)\right|\\
&\leq&
C\frac{kq}{n}\left[b_n^2\phi(h_n)v_n+b_n^2\phi(h_n)^{-1}h_n^{-2}v_n^{-a}\right].
\end{eqnarray*}
Then as in (\ref{eq4.12}), it follows that
\begin{eqnarray}\label{eq4.22}
k\sum_{1\leq i<j\leq q}|\mathrm{Cov}(Z_{ni},Z_{nj})|\rightarrow
0,\quad n\rightarrow \infty.
\end{eqnarray}
Once again, by stationarity,
\begin{eqnarray*}
\sum_{1\leq i<j\leq k} |\mathrm{Cov}(\xi_i,\xi_j)|&=&
\sum_{l=1}^{k-1} (k-l)|\mathrm{Cov}(\xi_1,\xi_{l+1})|\\
&\leq& k\sum_{l=1}^{k-1}|\mathrm{Cov}(\xi_1,\xi_{l+1})|\\
&\leq&
pk\sum_{l=1}^{k-1}\sum_{r=l(p+q)-q+1}^{l(p+q)+q-1}|\mathrm{Cov}(Z_{n1},Z_{n,r+1})|.
\end{eqnarray*}
Then, by Lemma \ref{lem4.1},
\begin{eqnarray*}
\sum_{1\leq i<j\leq k} |\mathrm{Cov}(\xi_i,\xi_j)|&\leq& C\frac{b_n^2 pk}{nh_n^2\phi(h_n)}\lambda_{p}\\
&\leq& C \log^2(n)\phi(h_n)^{-(1+\frac{2}{b})}\lambda_{p}.
\end{eqnarray*}
Now, define $p$ and $q$ as follows
\begin{eqnarray*}
p\sim \phi(h_n)^{-\delta_1},\quad q\sim
\phi(h_n)^{-\delta_2},\quad 0<\delta_2<\delta_1<\delta.
\end{eqnarray*}
We can choose $\delta_1$ and $\delta_2$ such that
$a>\frac{2+b}{\delta_2 b}$, then
\begin{eqnarray}\label{eq4.23}
\sum_{1\leq i<j\leq k} |\mathrm{Cov}(\xi_i,\xi_j)|\leq
C\log^2(n)\phi(h_n)^{a\delta_1-(1+\frac{2}{b})}\rightarrow 0,\quad
n\rightarrow \infty.
\end{eqnarray}
From (\ref{eq4.19}), (\ref{eq4.21}), (\ref{eq4.22}) and
(\ref{eq4.23}) it follows that
\begin{eqnarray*}
\mathrm{E}(T_n^{'})^2\rightarrow 0,\quad n\rightarrow \infty.
\end{eqnarray*}
\begin{eqnarray*}
\mathrm{E}[T_n^{''}]^2&\leq& (n-k(p+q))\mathrm{Var}(Z_{n1})+2
\sum_{1\leq i<j\leq
n} |\mathrm{Cov}(Z_{ni},Z_{nj})|\\
&\leq& p\mathrm{Var}(Z_{n1})+2 \sum_{1\leq i<j\leq
n} |\mathrm{Cov}(Z_{ni},Z_{nj})|\\
&\leq&
\frac{p\phi(h_n)}{n\left[\mathrm{E}\Delta_1(x)\right]^2}\mathrm{Var}(\hat{\Gamma}_1(x))+
\frac{C}{n\phi(h_n)}\sum_{1\leq i<j\leq
n}\mathrm{Cov}\left(\hat{\Gamma}_i(x),\hat{\Gamma}_j(x)\right).
\end{eqnarray*}
By a similar argument we find using (\ref{eq4.10}) and
(\ref{eq4.12}),
\begin{eqnarray*}
\mathrm{E}(T_n^{''})^2\rightarrow 0,\quad n\rightarrow \infty.
\end{eqnarray*}

\noindent \emph{Proof of convergence in (\ref{eq4.17})}.
\begin{align}\label{eq4.24}
\Big|\mathrm{E}\left(e^{it\sum_{j=1}^k
\eta_j}\right)-\prod_{j=1}^{k}\mathrm{E}\left(e^{it
\eta_j}\right)\Big|\leq&~\Big|\mathrm{E}\left(e^{it\sum_{j=1}^k
\eta_j}\right)-\mathrm{E}\left(e^{it\sum_{j=1}^{k-1}
\eta_j}\right)\mathrm{E}\left(e^{it \eta_k}\right)\Big|\nonumber\\
&+ \Big|\mathrm{E}\left(e^{it\sum_{j=1}^{k-1}
\eta_j}\right)-\prod_{j=1}^{k-1}\mathrm{E}\left(e^{it \eta_j}\right)\Big|\\
=&~\Big|\mathrm{Cov}\left(e^{it\sum_{j=1}^{k-1} \eta_j},e^{it
\eta_k}\right)\Big|+ \Big|\mathrm{E}\left(e^{it\sum_{j=1}^{k-1}
\eta_j}\right)-\prod_{j=1}^{k-1}\mathrm{E}\left(e^{it
\eta_j}\right)\Big|\nonumber.
\end{align}
By a repetition of this argument, inequality (\ref{eq4.24})
becomes
\begin{align}\label{eq4.25}
\Big|\mathrm{E}\left(e^{it\sum_{j=1}^k
\eta_j}\right)-\prod_{j=1}^{k}\mathrm{E}\left(e^{it
\eta_j}\right)\Big|\leq&~
\Big|\mathrm{Cov}\left(e^{it\sum_{j=1}^{k-1} \eta_j},e^{it
\eta_k}\right)\Big|+ \Big|\mathrm{Cov}\left(e^{it\sum_{j=1}^{k-2}
\eta_j},e^{it
\eta_{k-1}}\right)\Big|\nonumber\\
&+ \cdots+\Big|\mathrm{Cov}\left(e^{it \eta_2},e^{it
\eta_1}\right)\Big|.
\end{align}
Apply Lemma \ref{lem4.1} to each term on the right-hand side of
(\ref{eq4.25}) in order to obtain
\begin{align}\label{eq4.26}
\Big|\mathrm{E}\left(e^{it\sum_{j=1}^k
\eta_j}\right)-\prod_{j=1}^{k}\mathrm{E}\left(e^{it
\eta_j}\right)\Big|\leq&~Ct^2  \frac{\phi(h_n)b_n^2}{nh_n^2
\left[\mathrm{E}\Delta_1(x)\right]^2}\Big[\sum_{i\in
I_1}\sum_{j\in I_2} \lambda_{ij}+\sum_{i\in (I_1\cup
I_2)}\sum_{j\in I_3}
\lambda_{ij}\nonumber\\
&+ \cdots+\sum_{i\in (I_1\cup\ldots \cup I_{k-1})}\sum_{j\in I_k}
\lambda_{ij}\Big].
\end{align}
By stationarity, the inequality (\ref{eq4.26}) becomes
\begin{align}\label{eq4.27}
\Big|\mathrm{E}\left(e^{it\sum_{j=1}^k
\eta_j}\right)-\prod_{j=1}^{k}\mathrm{E}\left(e^{it
\eta_j}\right)\Big|\leq&~Ct^2\frac{\phi(h_n)b_n^2}{nh_n^2
\left[\mathrm{E}\Delta_1(x)\right]^2}\Big[(k-1)\sum_{i\in
I_1}\sum_{j\in I_2}
\lambda_{ij}\nonumber\\
&+ (k-2)\sum_{i\in I_1}\sum_{j\in I_3} \lambda_{ij}
+\cdots+\sum_{i\in I_1}\sum_{j\in I_k} \lambda_{ij}\Big].
\end{align}
Once again, by stationarity, for every $2\leq l\leq k$,
\begin{align*}
\sum_{i\in I_1}\sum_{j\in I_l} \lambda_{ij}=&~\lambda_{1,(l-1)(p+q)-p+2}+2\lambda_{1,(l-1)(p+q)-p+3}+\cdots+(p-1)\lambda_{1,(l-1)(p+q)}\\
&+ p~ \lambda_{1,(l-1)(p+q)+1}+(p-1)
\lambda_{1,(l-1)(p+q)+2}+\cdots+ \lambda_{1,(l-1)(p+q)+p}.
\end{align*}
Therefore, inequality (\ref{eq4.27}) becomes
\begin{eqnarray*}
\Big|\mathrm{E}\left(e^{it\sum_{j=1}^k
\eta_j}\right)-\prod_{j=1}^{k}\mathrm{E}\left(e^{it
\eta_j}\right)\Big|&\leq& Ct^2  \frac{\phi(h_n)b_n^2}{nh_n^2
\left[\mathrm{E}\Delta_1(x)\right]^2}pk \lambda_{q}\nonumber\\
&\leq&Ct^2  \frac{b_n^2 pk}{nh_n^2 \phi(h_n)}\lambda_{q}\\
&\leq& C\log^2(n)\phi(h_n)^{a\delta_2-(1+\frac{2}{b})}\rightarrow
0,\quad n\rightarrow \infty.
\end{eqnarray*}

\noindent \emph{Proof of convergence in (\ref{eq4.18})}. By (\ref{eq4.10}) and (\ref{eq4.22}), we have
\begin{eqnarray*}
k\mathrm{Var}(\eta_1)=kp\mathrm{Var}(Z_{n1})+2k\sum_{1\leq i<j\leq
p} |\mathrm{Cov}(Z_{ni},Z_{nj})|\rightarrow \sigma_1^2(x).
\end{eqnarray*}
Then, since $\ds|\eta_1|\leq C\frac{b_n p}{ (n\phi(h_n))^{1/2}}$,
it follows that
\begin{eqnarray*}
k\mathrm{E}\left(\eta_1^{2}\mathds{1}_{\{\eta_1>\varepsilon
\sigma_1(x)
\}}\right)&\leq& C\frac{kb_n^2p^2}{n\phi(h_n)}P(\eta_1>\varepsilon \sigma_1(x))\\
&\leq& C \frac{b_n^2p^2}{n\phi(h_n)}\frac{k\mathrm{Var(\eta_1)}}{\varepsilon^2\sigma_1^2(x)}\\
&\leq&
C\frac{\log^2(n)}{n\phi(h_n)^{1+2\delta}}\frac{k\mathrm{Var(\eta_1)}}{\varepsilon^2\sigma_1^2(x)}\rightarrow
0.
\end{eqnarray*}
The second assertion is an application of the first one when
replacing $\phi(Y_i)$ with $\phi(Y_i)-r(x)$.\cqfd \vspace{0.5cm}

\noindent\textbf{Proof of Corollaire \ref{thm4.2}.} Consider the following decomposition:
\begin{align*}
r_n(x)-r(x)=&~(g_n-rf_n)(x)-\mathrm{E}((g_n-rf_n)(x))\\
&- (r(x)-\mathrm{E}g_n(x)) + (r(x)-\mathrm{E}g_n(x))(f_n(x)-\mathrm{E}f_n(x))\\
&-
(g_n(x)-\mathrm{E}g_n(x))(f_n(x)-\mathrm{E}f_n(x))+r_n(x)(f_n(x)-\mathrm{E}f_n(x))^2.
\end{align*}
Since
\begin{eqnarray*}
(n\phi(h_n))^{1/2}\left[r(x)-\mathrm{E}g_n(x)\right]=\mathcal{O}\left(\left(n\phi(h_n)
h_n^{2\beta}\right)^{1/2}\right),
\end{eqnarray*}
\begin{eqnarray*}
(n\phi(h_n))^{1/2}\left(r(x)-\mathrm{E}g_n(x)\right)\mathrm{E}\left[(f_n(x)-\mathrm{E}f_n(x))\right]=\mathcal{O}\left(h_n^{\beta}\right),
\end{eqnarray*}
\begin{eqnarray*}
\lefteqn{(n\phi(h_n))^{1/2}\mathrm{E}\left|(g_n(x)-\mathrm{E}g_n(x))(f_n(x)-\mathrm{E}f_n(x))\right|}\\
&\leq&(n\phi(h_n))^{1/2}\left[\mathrm{E}(g_n(x)-\mathrm{E}g_n(x))^2\right]^{1/2}\left[\mathrm{E}(f_n(x)-\mathrm{E}f_n(x))^2\right]^{1/2}\\
&=& \mathcal{O}\left(\left(n\phi(h_n)\right)^{-1/2}\right)
\end{eqnarray*}
and
\begin{eqnarray*}
(n\phi(h_n))^{1/4}\mathrm{E}\left|r_n(x)(f_n(x)-\mathrm{E}f_n(x))^2\right|^{1/2}&\leq&
\left(n\phi(h_n)\right)^{1/4}\left(\mathrm{E}|r_n(x)|\right)^{1/2}\left(\mathrm{E}(f_n(x)-\mathrm{E}f_n(x))^2\right)^{1/2}\\
&=&\mathcal{O}\left(\log^{1/2}(n)\left(n\phi(h_n)\right)^{-1/4}\right).
\end{eqnarray*}
The result then follows from Theorem \ref{prop4.3}.\cqfd
\newpage

\bibliographystyle{natbib}
\bibliography{reg}

\end{document}